\documentclass[11pt]{amsart}
\usepackage{amssymb,color}
\usepackage{graphicx}

\numberwithin{equation}{section}

\newcommand{\bea}{\begin{eqnarray}}
\newcommand{\eea}{\end{eqnarray}}
\newcommand{\ba}{\begin{array}}
\newcommand{\ea}{\end{array}}
\newcommand{\edc}{\end{document}}
\newcommand{\bc}{\begin{center}}
\newcommand{\ec}{\end{center}}
\newcommand{\be}{\begin{equation}}
\newcommand{\ee}{\end{equation}}

\def\bc{{\mathbb C}}

\newtheorem{thm}{Theorem}[section]

\newtheorem{defin}[thm]{Definition}
\theoremstyle{remark}
\newtheorem{rem}{Remark}[section]
\newtheorem{ex}{Example}[section]

\date{\today}%
\begin{document}

\title[Cellular automata on Cayley tree]
{The entropy and reversibility of cellular automata on Cayley
tree}

\author{Hasan Ak\i n}
\address{Hasan Ak\i n, Ceyhun
Atuf Kansu Caddesi 1164. Sokak, 9/4, TR06105, \c{C}ankaya, Ankara,
Turkey} \email{{\tt akinhasan25@gmail.com}}
\begin{abstract}
In this paper, we study linear cellular automata (CAs) on Cayley
tree of order 2 over the field $\mathbb F_p$ (the set of prime
numbers modulo $p$). We construct the rule matrix corresponding to
finite cellular automata on Cayley tree. Further, we analyze the
reversibility problem of this cellular automata for some given
values of $a,b,c,d\in \mathbb{F}_{p}\setminus \{0\}$ and the
levels $n$ of Cayley tree. We compute the measure-theoretical
entropy of the cellular automata which we define on Cayley tree.
We show that for CAs on Cayley tree the measure entropy with
respect to uniform Bernoulli measure is infinity.

\end{abstract}

\maketitle 
\section{Introduction}
A cellular automaton (plural cellular automata, shortly CA) has
been studied and applied as a discrete model in many areas of
science. Cellular automata (CAs) have very rich computational
properties and provide different models in computation. CAs were
first used for modeling various physical and biological processes
and especially in computer science. Recently,  CAs have been
widely investigated in many disciplines with different purposes
such as simulation of natural phenomena, pseudo-random number
generation, image processing, analysis of a universal model of
computations, coding theory, cryptography, ergodic theory
(\cite{A1,CCD,KCDMS,KCDV,YZP,DCK}).

Most of the studies and applications for CA is extensively done
for one-dimensional (1-D) CA. "The Game of Life" developed by John
H. Conway in the 1960Тs is an example of a two-dimensional (2-D)
CA. John von Neumann in the late 40's and early 50's studied CA as
a self-reproducing simple organisms \cite{John}. 2-D CA with von
Neumann neighborhood has found many applications and been explored
in the literature \cite{SOMA}. Nowadays, 2-D CAs have attracted
much of the interest. Some basic and precise mathematical models
using matrix algebra built on field $\mathbb Z_2$  were reported
for characterizing the behavior of two-dimensional nearest
neighborhood linear CAs with null or periodic boundary conditions
\cite{CCD,KCDMS,KCDV,DCK,SOMA}.

The reversibility problem of some special classes of 1-D CAs
reflective and periodic boundary conditions has been studied with
the help of matrix algebra approach by several researchers
\cite{HFI,siap}. In \cite{ASU2014a}, we have defined a family of
one-dimensional finite linear cellular automata with reflective
boundary condition over the field $\mathbb Z_p$. In \cite{SUAS},
we investigated 2D finite CA with a von Neumann neighborhood under
periodic, adiabatic or reflexive boundaries conditions over the
ternary field the field $\mathbb Z_3$, which can be considered as
a three-state case. The application of linear rules on image
matrix is demonstrated which forms the basis of self-replicating
and self-similar patterns in image processing
\cite{SUAS,USAS,USAS1,SUA}. Particularly the rules are used for
image multiplication of one image into several replicating or
similar images. In \cite{KSA}, we investigated error correcting
codes via reversible cellular automata over finite fields. In this
paper, we start with linear cellular automata (CA) in relation to
a basic mathematical structure on regular Cayley tree of order 2.
Recently, we have investigated the reversibility problem of
multidimensional linear cellular automata under certain boundary
conditions (null, periodic, reflective) on some lattices; however,
we have not obtained exact algorithms for determining whether a
multidimensional linear cellular automaton is reversible
\cite{ASU2014a,HFI,siap,ChangAkin2017,KSA,ASSY,KorogluSA2016}.

In Ref. \cite{FF}, Fici and Fiorenzi have a first attempt to study
topological properties of CA on the full tree shift $A^{\sum
^{*}}$, where $\sum ^{*}$ is the free monoid of finite rank
$|\sum|$. In this case, the Cayley graph of $\sum ^{*}$ is a
regular $|\sum|$-ary rooted tree. Fici and Fiorenzi \cite{FF} have
studied cellular automata defined on the full $k$-ary tree shift
(for $k \geq 2$). In this paper, we study cellular automata on
regular Cayley tree of order 2.

Several notions of the entropy of measure-preserving
transformation on probability space in ergodic theory have been
investigated \cite{AMalay,A4}. The notion of entropy, both
topological and measure-theoretical is one of the fundamental
invariants in ergodic theory. In the last years, a lot of works
have been devoted to this subject \cite{A1,AMalay,A2,A3,DMM}.
Recall that by the Variational Principle the topological entropy
is the supremum of the entropies of invariant measures. In
\cite{A1}, the author has shown that the uniform Bernoulli measure
is a measure of maximal entropy for some 1-D LCAs. Morris and Ward
\cite{MW} proved that an ergodic additive CA in two dimensions has
infinite topological entropy (see \cite{Mey} for details).
Recently, Blanchard and Tisseur \cite{BT} have introduced the
entropy rate of multidimensional CAs and proved several results
that show that entropy rate of 2-D CA preserve similar properties
of the entropy of 1-D CA.

In this present paper, firstly we define cellular automata on
Cayley tree (or Bethe lattice) of order 2. This generalizes the
case of one-sided CA (where order of the Cayley tree is one). We
construct a transition rule matrix corresponding to finite
cellular automata on Cayley tree by using matrix algebra built on
the field $\mathbb Z_p$ (the set of prime numbers modulo $p$).
Further, we discuss the reversibility problem of this cellular
automata. Lastly, we study the measure theoretical entropy of the
CAs on Cayley tree. We show that for CAs on Cayley tree the
measure entropy with respect to uniform Bernoulli measure is
infinity.
\section{Finite CA over Cayley tree}
\begin{figure} [!htbp]%
\centering
\includegraphics[width=45mm]{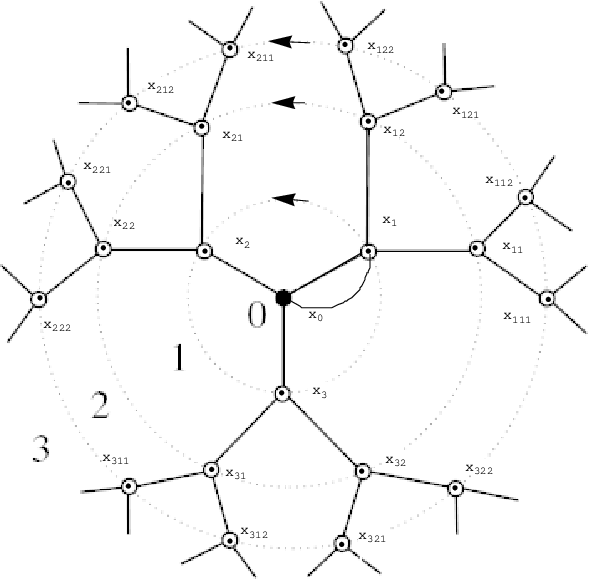}\ \ \ \ \ \ \ \ \ \ \ \
\includegraphics[width=40mm]{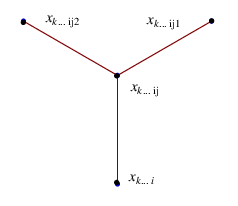}
\caption{a) Cayley tree of order two with levels 3, b) Elements of
the nearest neighborhoods surround the center $x_{k\ldots ij}$,
$k=1,2,3$ and $i,j=1,2$. }\label{fig2}
\end{figure}

Let $\mathbb{F}_{p}=\{0, 1,\ldots, p-1\}$ $(p\geq 2)$ be the field
of the prime numbers modulo $p$ ($\mathbb{F}_{p}$ is called a
\emph{state space}). The Cayley tree $\Gamma^{k}$ of order
$k\geq1$ is an infinite tree, i.e., a graph without cycles, from
each vertex of which exactly $k+1$ edges issue. Let
$\Gamma^{k}=(V, L, i)$, where $V$ is the set of vertices of
$\Gamma^{k}$, $L$ is the set of edges of $\Gamma^{k}$ and $i$ is
the incidence function associating each edge $\ell\in L$ with its
end points $x, y \in V.$ A configuration $\sigma$ on $V$ is
defined as a function $x\in V\to\sigma (x)\in\mathbb F_p $; in a
similar manner one defines configurations $\sigma_n$ and $\omega$
on $V_n$ and $W_n$, respectively. The set of all configurations on
$V$ (resp. $V_n$, $W_n$) coincides with $\Omega=\mathbb F_p^{V}$
(resp. $\Omega_{V_n}=\mathbb F_p^{V_n},\ \ \Omega_{W_n}=\mathbb
F_p^{W_n}$). One can see that
$\Omega_{V_n}=\Omega_{V_{n-1}}\times\Omega_{W_n}$. Denote by
$\mathbb{F}^{\Gamma^{2}}_{p}$, i.e., the set of all configurations
on $\Gamma^{2}$. In the sequel we will consider Cayley tree
$\Gamma^{2}=(V, L,i)$ with the root $x_0$. If $i(\ell)=\{x, y\}$,
then $x$ and $y$ are called the nearest neighboring vertices and
we write $\ell=<x, y>$. For $x,\ y \in V$, the distance $d(x, y)$
on Cayley tree is defined by the formula
\begin{eqnarray*}
d(x, y)&=&\min \{d|x=x_{0}, x_{1}, x_{2}, \ldots, x_{d-1}, x_{d}=y
\in V
\ \text{such  that  the pairs} \\
&&  <x_0, x_1>, \ldots, <x_{d-1} , x_d> \text{are neighboring
vertices}\}.
\end{eqnarray*}
For the fixed root vertex $x^{0}\in V$ we have
$$
W_{n}={\{}x\in V: d(x^{0},x)=n{\}}
$$
$$
V_{n}=\{x\in V: d(x^{0},x)\leq n\},
$$
$$
L_{n}=\{\ell= <x,y>\in L:x,y\in V_{n} \}.
$$
In this section, we will order the elements of $V_{n}$ in the
lexicographical meaning (see \cite{LMS}) as the Fig. \ref{fig2}.
Given two vertices $x_{u},x_{v}$, the lexicographical order of
$x_{u},x_{v}$ is defined as $x_{u}\preceq x_{v}$ if and only if $u
\preceq v$.

Let us rewrite the elements of $W_{n}$ in the following order,
$$
\overrightarrow{W_{n}}:=(x^{(1)}_{W_{n}},x^{(2)}_{W_{n}},\ldots,x^{(|W_{n}|)}_{W_{n}}).
$$
One can easily compute equations $|W_{n}|=3.2^{(n-1)}$ and
$|V_{n}|=1+3.(2^{n}-1)$. For the sake of shortness, throughout the
paper we are going to represent vertices
$x^{(1)}_{W_{n}},x^{(2)}_{W_{n}},\ldots,x^{(|W_{n}|)}_{W_{n}}$ of
$W_{n}$ by means of the coordinate system as follows:
\begin{eqnarray*}
&&x^{(1)}_{W_{n}}=x_{11...11},x^{(2)}_{W_{n}}=x_{11...12},x^{(3)}_{W_{n}}=x_{11...21},x^{(4)}_{W_{n}}=x_{11...22},\\
&&...\\
&&x^{(|W_{n}|-3)}_{W_{n}}=x_{32...211},x^{((|W_{n}|-2))}_{W_{n}}=x_{32...212},
x^{(|W_{n}|-1))}_{W_{n}}=x_{32...221},x^{(|W_{n}|))}_{W_{n}}=x_{32...222}.
\end{eqnarray*}
In the Fig. \ref{fig2}, we show Cayley tree of order two with
levels 3 and the nearest neighborhood which comprises three cells
which surround the center cell $x_{k\ldots ij}$. The state
$x_{k\ldots ij}^{(t+1)}$ of the cell $(i,j)^{\text{th}}$ at time
$(t+1)$ is defined by the local rule function $f:\mathbb{F}_p^4
\rightarrow \mathbb{F}_p$ as follows:
\begin{eqnarray}\label{CA}
x_{k\ldots ij}^{(t+1)} &=&f(x_{k\ldots i}^{(t)},x_{k\ldots
ij}^{(t)},x_{k\ldots ij1}^{(t)},x_{k\ldots ij2}^{(t)})\\\nonumber
&=& ax_{k\ldots ij1}^{(t)}+bx_{k\ldots ij2}^{(t)}+ cx_{k\ldots
i}^{(t)}+dx_{k\ldots ij}^{(t)} (\text{mod }p),
\end{eqnarray}
where $a,b,c,d\in \mathbb{F}_{p}\setminus \{0\}$, $x_{k\ldots
i}^{(t)}\in W_{n-2}$, $x_{k\ldots ij}^{(t)}\in W_{n-1}$ and
$x_{k\ldots ij1}^{(t)}, x_{k\ldots ij2}^{(t)}\in W_{n}$, $k=1,2,3$
and $i,j=1,2$ (see the Fig. \ref{fig2} (b)).\\
Specifically, for state $x_{0}^{(t+1)}$ in the root vertex we can
show
\begin{eqnarray}\label{CA1}
x_{0}^{(t+1)}
&=&f(x_{0}^{(t)},x_{1}^{(t)},x_{2}^{(t)},x_{3}^{(t)})=
ax_{1}^{(t)}+bx_{2}^{(t)}+ cx_{3}^{(t)}+dx_{0}^{(t)} (\text{mod
}p).
\end{eqnarray}
Function
\begin{eqnarray}\label{mapI}
T_f:\mathbb{F}^{\Gamma^{2}}_{p}\rightarrow
\mathbb{F}^{\Gamma^{2}}_{p}
\end{eqnarray}
is called a cellular automaton (CA) generated by the rules
\eqref{CA} and \eqref{CA1}.
If the boundary cells are connected to 0-state, then CA are called
Null Boundary CA, i.e., $V\setminus W_{n}=\{0\}$ for a fixed $n$.
If the same rule is applied to all of the cells in ever
evaluation, then those CA are called uniform or regular.

In Sections \ref{rule matrix} and \ref{Reversibility}, we consider
linear transformations of finite dimensional vectors spaces
corresponding to these finite linear cellular automata by imposing
the null boundary condition, which means that the states of cells
outside a given ball around the origin are fixed to be zero.

\section{Construction of the rule matrix in the finite
case}\label{rule matrix} In this section, we can characterize
finite cellular automata with Null boundary condition over Cayley
tree of order two over the field $\mathbb{F}_{p}$. In order to
characterize the corresponding rule, first we represent finite
Cayley tree $n$ level as a column vector of size
$(1+3(2^{n}-1))\times 1.$

Let us denote all configurations of Cayley tree with levels $n$ by
$\Omega_n$. In order to accomplish this goal we define the
following map
$$
\Phi:\Omega_n\to \mathbb{M}_{(1+3(2^{n}-1))\times 1}(\mathbb Z_p),
$$
which takes  the $t^{th}$ state $X^{(t)}$ given by
$$\Omega_n
\rightarrow X^{(t)}:=\left(x_{0}^{(t)},x_{1}^{(t)},\ldots
,x_{21...11}^{(t)},x_{21...12}^{(t)}, \ldots ,x_{32...211}^{(t)},
 x_{32...212}^{(t)},x_{32...221}^{(t)},x_{32...222}^{(t)}\right)^T,
$$
where the superscript $T$ denotes the transpose and
$\mathbb{M}_{(1+3(2^{n}-1))\times 1}(\mathbb Z_p)$ is the set of
matrices with  entries $\mathbb F_p$.

The configuration $\sigma^{(t)}_n\in \Omega_n$ is called the
configuration matrix (or information matrix) of the finite CA on
Cayley tree with levels $n$ at time $t$ and $\sigma^{(0)}_n$ is
initial information matrix of the finite CA. The whole evolution
of a particular cellular automata can be comprised in its global
transition function \cite{DCK} (see \cite{DCK,SAS1,IAS} for the
square lattice $\mathbb F^2$ and see \cite{SAS} for the hexagonal
lattice).

Therefore, one can conclude that $\Phi(\sigma^{(t)}_n) =
X^{(t)}_{(1+3(2^{n}-1))\times 1}$. Using the identification
\eqref{mapI}, due to linearity of the finite CA we can define as
follows:
$$
(M_{R}^{(n)})_{(1+3(2^{n}-1))\times
(1+3(2^{n}-1))}X^{(t)}_{(1+3(2^{n}-1))\times
1}=X^{(t+1)}_{(1+3(2^{n}-1))\times 1},
$$
where $n$ is the number of levels of the Cayley tree.
\begin{thm}\label{main}
Let $a,b,c,d\in  \mathbb{F}_p^{*}=\mathbb{F}_p \setminus \{0 \},$
$n\geq 2$. Then, the transition rule matrix\\
$(M_{R}^{(n)})_{1+3(2^{n}-1)\times 1+3(2^{n}-1)}$ corresponding to
the finite cellular automata on Cayley tree of order two with
$n$-level finite over NB is given by 
{\small \footnotesize \begin{equation}\label{rm}
\left(
\begin{array}{cccccccccccccccc}
  d & P & 0_{1\times6} & \cdots & 0 & 0 & 0 \\
  Q & D_{3\times3} &B_{3\times6}&  0  & \cdots & 0 & 0 \\
  0 & C_{6\times3}& D_{6\times6} &  B_{6\times12}& 0 & \cdots  & 0   \\
  \vdots & \vdots & \vdots & \ddots & \ddots & \ddots & \vdots \\
 0& \cdots & 0 & C_{3.2^{n-3}\times3.2^{n-4}}& D_{3.2^{n-3}\times3.2^{n-3}} &  B_{3.2^{n-3}\times3.2^{n-2}}& 0   \\
 0&0 &  \cdots   & 0 & C_{3.2^{n-2}\times3.2^{n-3}}& D_{3.2^{n-2}\times3.2^{n-2}} &  B_{3.2^{n-2}\times3.2^{n-1}}\\
 0 & 0 & 0  & \cdots & 0 & C_{3.2^{n-1}\times3.2^{n-2}}& D_{3.2^{n-1}\times3.2^{n-1}} \\
\end{array}
\right)
\end{equation}}
where each submatrices are as follows: $P=(\begin{array}{ccc}a& b
& c\end{array})$, $Q=\left(\begin{array}{ccc} c\\c\\c
\end{array}\right),$
{\small \footnotesize \begin{equation*}
C_{3.2^{n-i}\times 3.2^{n-(i+1)}}=\left(%
\begin{array}{cccccccc}
  c& 0 & 0 & 0  & \cdots & 0 & 0 \\
  c & 0 & 0& 0   & \cdots & 0 & 0 \\
  0 & c & 0 & 0 & \cdots & 0 & 0   \\
  0 & c & 0 & 0 & \cdots & 0 & 0   \\
\cdots & \cdots & \cdots & \cdots & \cdots & \cdots & \cdots \\
  0 & \cdots  &\cdots & 0 & 0 & c& 0 \\
  0 & 0&\cdots  &\cdots & 0 & c & 0\\
   0 & \cdots  &\cdots & 0 & 0 & 0& c \\
  0 & 0&\cdots  &\cdots & 0 & 0 & c\\
\end{array}%
\right),
\end{equation*}}
{\small \footnotesize  \begin{equation*}
B_{3.2^{n-(i+1)}\times3.2^{n-i}}=\left(%
\begin{array}{cccccccccc}
  a& b & 0 & 0  & 0& 0 & 0& \cdots & 0 \\
  0 & 0 & a& b   & 0& 0 & 0& \cdots & 0 \\
  0 & 0& 0 & 0  & a & b & 0& \cdots & 0  \\
\cdots & \cdots & \cdots & \cdots & \cdots & \cdots & \cdots \\
0 & \cdots  &\cdots & 0 & a & b& 0& 0 & 0 \\
  0 & 0&\cdots  &\cdots & 0 & 0 & 0& a & b\\
\end{array}%
\right),
\end{equation*}}
and
\begin{equation*} \label{pmnb}
D_{3.2^{n-i}\times3.2^{n-i}}=d.I_{3.2^{n-i}\times3.2^{n-i}}.
\end{equation*} For  $i=1,2,\ldots,n-1.$
 \end{thm}
In the Theorem \ref{main}, we have obtained a general form of the
matrix representation (or transition rule matrix) for these linear
transformations with respect to a basis given by the
lexicographical order on the vertices. We do not include a
detailed proof of the theorem which gives the rule matrix of CA.
The proof is obtained by determining the image of the basis
elements of the space $\mathbb F_p^{(1+3(2^{n}-1))}$ under the CA.
These images contribute to the columns of the rule matrix.

Let us illustrate this in Examples \ref{exRule} and \ref{exam2}.
\begin{ex}\label{exRule} If we take the number of level as $n=2$, then we get the rule
matrix $M_{R}$ of order 10. We consider a configuration
$\sigma^{(t)}_2$ of number of levels 2 with null boundary:
\begin{figure} [!htbp]%
\centering
\includegraphics[width=50mm]{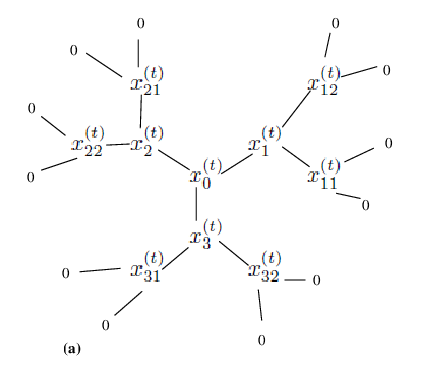}\ \ \ \
\includegraphics[width=50mm]{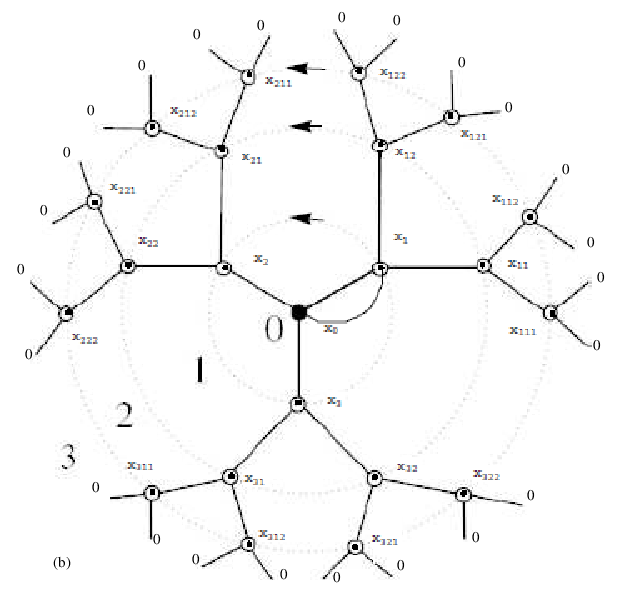}
\caption{A configuration $\sigma^{(t)}_2$ of levels 2 and 3 with
null boundary on Cayley tree of order two.}\label{level2}
\end{figure}
Let us apply the local rules \eqref{CA} and \eqref{CA1} on
configuration $\sigma^{(t)}_2$ in the Fig \ref{level2}. Then, we
get a new configuration under this transformation which is
\begin{eqnarray*}
x_{0}^{(t+1)}&=&ax_{1}^{(t)}+bx_{2}^{(t)}+
cx_{3}^{(t)}+dx_{0}^{(t)};\ x_{1}^{(t+1)}=ax_{11}^{(t)}+bx_{12}^{(t)}+cx_{0}^{(t)}+dx_{1}^{(t)}\\
x_{2}^{(t+1)}&=&ax_{21}^{(t)}+bx_{22}^{(t)}+cx_{0}^{(t)}+dx_{2}^{(t)};\ x_{3}^{(t+1)}=ax_{31}^{(t)}+bx_{32}^{(t)}+cx_{0}^{(t)}+dx_{3}^{(t)}\\
x_{11}^{(t+1)}&=&cx_{1}^{(t)}+dx_{11}^{(t)};\ x_{12}^{(t+1)}=cx_{1}^{(t)}+dx_{12}^{(t)};\ x_{21}^{(t+1)}=cx_{2}^{(t)}+dx_{21}^{(t)}\\
x_{22}^{(t+1)}&=&cx_{2}^{(t)}+dx_{22}^{(t)};\
x_{31}^{(t+1)}=cx_{3}^{(t)}+dx_{31}^{(t)};\
x_{32}^{(t+1)}=cx_{3}^{(t)}+dx_{32}^{(t)}.
\end{eqnarray*}
Hence, we obtain the rule matrix $M_{R}^{(2)}$ of order 10 as
follows: {\small \footnotesize  \begin{eqnarray*}
M_{R}^{(2)}=\left(
\begin{array}{c|ccc|cccccc}
 d & a & b & c & 0 & 0 & 0 & 0 & 0 & 0 \\\hline
 c & d & 0 & 0 & a & b & 0 & 0 & 0 & 0 \\
 c & 0 & d & 0 & 0 & 0 & a & b & 0 & 0 \\
 c & 0 & 0 & d & 0 & 0 & 0 & 0 & a & b \\\hline
 0 & c & 0 & 0 & d & 0 & 0 & 0 & 0 & 0 \\
 0 & c & 0 & 0 & 0 & d & 0 & 0 & 0 & 0 \\
 0 & 0 & c & 0 & 0 & 0 & d & 0 & 0 & 0 \\
 0 & 0 & c & 0 & 0 & 0 & 0 & d & 0 & 0 \\
 0 & 0 & 0 & c & 0 & 0 & 0 & 0 & d & 0 \\
 0 & 0 & 0 & c & 0 & 0 & 0 & 0 & 0 & d
\end{array}
\right).
\end{eqnarray*}}
\end{ex}
\begin{ex}\label{exam2} Let us consider the configuration with levels 3 given in the Fig.
\ref{level2}.
If we apply the rules \eqref{CA} and \eqref{CA1}, then we obtain
the following rule matrix $M_R^{(3)}$: {\small \footnotesize
\begin{eqnarray*}\label{matn=3}
M_R^{(3)}&=&\left(
\begin{array}{cccccccccccccccc}
 d & P & 0_{1\times6} & 0 \\
 Q & D_{3\times3} &B_{3\times6}&  0 \\
 0 & C_{6\times3}& D_{6\times6} &  B_{6\times12}\\
 0 &0 & C_{12\times6}& D_{12\times12} \\
\end{array}
\right)\\\nonumber
 &=&
\left(
\begin{array}{c|ccc|cccccc|cccccccccccc}
 d & a & b & c & 0 & 0 & 0 & 0 & 0 & 0 & 0 & 0 & 0 & 0 & 0 & 0 & 0 & 0 & 0 & 0 & 0 & 0 \\\hline
 c & d & 0 & 0 & a & b & 0 & 0 & 0 & 0 & 0 & 0 & 0 & 0 & 0 & 0 & 0 & 0 & 0 & 0 & 0 & 0 \\
 c & 0 & d & 0 & 0 & 0 & a & b & 0 & 0 & 0 & 0 & 0 & 0 & 0 & 0 & 0 & 0 & 0 & 0 & 0 & 0 \\
 c & 0 & 0 & d & 0 & 0 & 0 & 0 & a & b & 0 & 0 & 0 & 0 & 0 & 0 & 0 & 0 & 0 & 0 & 0 & 0 \\\hline
 0 & c & 0 & 0 & d & 0 & 0 & 0 & 0 & 0 & a & b & 0 & 0 & 0 & 0 & 0 & 0 & 0 & 0 & 0 & 0 \\
 0 & c & 0 & 0 & 0 & d & 0 & 0 & 0 & 0 & 0 & 0 & a & b & 0 & 0 & 0 & 0 & 0 & 0 & 0 & 0 \\
 0 & 0 & c & 0 & 0 & 0 & d & 0 & 0 & 0 & 0 & 0 & 0 & 0 & a & b & 0 & 0 & 0 & 0 & 0 & 0 \\
 0 & 0 & c & 0 & 0 & 0 & 0 & d & 0 & 0 & 0 & 0 & 0 & 0 & 0 & 0 & a & b & 0 & 0 & 0 & 0 \\
 0 & 0 & 0 & c & 0 & 0 & 0 & 0 & d & 0 & 0 & 0 & 0 & 0 & 0 & 0 & 0 & 0 & a & b & 0 & 0 \\
 0 & 0 & 0 & c & 0 & 0 & 0 & 0 & 0 & d & 0 & 0 & 0 & 0 & 0 & 0 & 0 & 0 & 0 & 0 & a & b \\\hline
 0 & 0 & 0 & 0 & c & 0 & 0 & 0 & 0 & 0 & d & 0 & 0 & 0 & 0 & 0 & 0 & 0 & 0 & 0 & 0 & 0 \\
 0 & 0 & 0 & 0 & c & 0 & 0 & 0 & 0 & 0 & 0 & d & 0 & 0 & 0 & 0 & 0 & 0 & 0 & 0 & 0 & 0 \\
 0 & 0 & 0 & 0 & 0 & c & 0 & 0 & 0 & 0 & 0 & 0 & d & 0 & 0 & 0 & 0 & 0 & 0 & 0 & 0 & 0 \\
 0 & 0 & 0 & 0 & 0 & c & 0 & 0 & 0 & 0 & 0 & 0 & 0 & d & 0 & 0 & 0 & 0 & 0 & 0 & 0 & 0 \\
 0 & 0 & 0 & 0 & 0 & 0 & c & 0 & 0 & 0 & 0 & 0 & 0 & 0 & d & 0 & 0 & 0 & 0 & 0 & 0 & 0 \\
 0 & 0 & 0 & 0 & 0 & 0 & c & 0 & 0 & 0 & 0 & 0 & 0 & 0 & 0 & d & 0 & 0 & 0 & 0 & 0 & 0 \\
 0 & 0 & 0 & 0 & 0 & 0 & 0 & c & 0 & 0 & 0 & 0 & 0 & 0 & 0 & 0 & d & 0 & 0 & 0 & 0 & 0 \\
 0 & 0 & 0 & 0 & 0 & 0 & 0 & c & 0 & 0 & 0 & 0 & 0 & 0 & 0 & 0 & 0 & d & 0 & 0 & 0 & 0 \\
 0 & 0 & 0 & 0 & 0 & 0 & 0 & 0 & c & 0 & 0 & 0 & 0 & 0 & 0 & 0 & 0 & 0 & d & 0 & 0 & 0 \\
 0 & 0 & 0 & 0 & 0 & 0 & 0 & 0 & c & 0 & 0 & 0 & 0 & 0 & 0 & 0 & 0 & 0 & 0 & d & 0 & 0 \\
 0 & 0 & 0 & 0 & 0 & 0 & 0 & 0 & 0 & c & 0 & 0 & 0 & 0 & 0 & 0 & 0 & 0 & 0 & 0 & d & 0 \\
 0 & 0 & 0 & 0 & 0 & 0 & 0 & 0 & 0 & c & 0 & 0 & 0 & 0 & 0 & 0 & 0 & 0 & 0 & 0 & 0 & d
\end{array}
\right).
\end{eqnarray*}}
\end{ex}
In order to illustrate the behavior of the finite CA on Cayley
tree, we can study the image and preimage under finite CA of a
configuration by means of relating matrix and its inverse matrix
(see \cite{FF}).

\section{Reversibility of CA on Cayley tree with Null
Boundary}\label{Reversibility} In this section, we characterize
finite cellular automata with NBC determined by nearest neighbor
rule on Cayley tree. For finite CA, in order to obtain the
reversible of a finite CA many authors
\cite{KCDMS,KCDV,YZP,HFI,siap,Rey,IAS,SAS} have used the rule
matrices. It is well known that a cellular automaton is reversible
if and only if it is bijective \cite{Kari}. Since we already have
found the rule matrix $M_{R}^{(n)}$ corresponding to the the
finite CA, by using the matrix in \eqref{rm}, we can state the
following relation between the column vectors $X^{(t)}$ and the
rule matrix $M_{R}$:
\[
X^{(t+1)}=M_{R}^{(n)}X^{(t)}\ (\text{mod}\ p).
\]
If the rule matrix $M_{R}^{(n)}$ is non-singular, then we have
\[
X^{(t)}=(M_{R}^{(n)})^{-1}X^{(t+1)}\ (\text{mod}\ p).
\]
Thus, in this paper, one of our main aims is to study whether the
rule matrix $M_{R}^{(n)}$ in \eqref{rm} is invertible or not. It
is well known that the finite CA is reversible if and only if its
rule matrix $M_{R}^{(n)}$ is non-singular (see
\cite{KCDMS,KCDV,YZP,IAS,SAS} for details). If the determinant of
a matrix is not equal to zero, then it is invertible, so the CA on
Cayley tree is reversible, otherwise, it is irreversible. If the
CA is not invertible, then one can study "Garden of Eden" for the
finite CA (see \cite{FF,SAS1}).

It is well known that the 1D finite CA is reversible iff its rule
matrix $M_R^{(n)}$ is non- singular (see \cite{siap} for details).
An efficient tool to compute the determinant of a matrix $A$ is to
multiply all eigenvalues of $A$. So, we conclude that the
reversibility of the original system comes from the combination of
eigenvalues of these components, so does its inverse
\cite{ChangAkin2017}.

Let us consider matrix $M_R^{(n)}$. The characteristic polynomial
of the matrix $M_R^{(n)}$ is given by

$$
\Delta _{M_R^{(n)}}(\lambda )=\det (\lambda I-M_R^{(n)})=\lambda
^n+\sum _{i=1}^n a_i\lambda ^{n-i}=\prod _{i=1}^n (\lambda
-\lambda _i).
$$
If we assume $\lambda=0$, then from the last equation we have
$$
\Delta _{M_R^{(n)}}(0)=\det (-M_R^{(n)})=(-1)^{n}\det
(M_R^{(n)})=(-1)^n\prod _{i=1}^n \lambda _i=a_n.
$$
Therefore, if $\det(M_R^{(n)})\neq 0$, then $M_R^{(n)}$ is
invertible, so corresponding CA is reversible. On the other hand,
due to $\det(M_R^{(n)})=\prod _{i=1}^n \lambda _i$, if 0 is not an
eigenvalue of $M_R^{(n)}$ over $\mathbb{F}_p$, then corresponding
CA is reversible, where for $(i=1,2,\cdots ,3.2^{n}-2)$ $\lambda
_i$ is an eigenvalue of $M_R^{(n)}$.

The following theorem provides basic transitions for reversibility
of 1D finite CA on Cayley tree of order 2.
\begin{thm}\label{theorem-hybrid}
The linear cellular automaton $T_{f}$ over $\mathbb{F}_p$ under
null boundary condition is characterized by the matrix $T_n$, and
vice versa. More explicitly, the diagram
\begin{eqnarray*}
\begin{array}{ccc}
 \mathbb{F}^{\Gamma^{2}}_{p} & \overset{T_{f}}{\longrightarrow } & \mathbb{F}^{\Gamma^{2}}_{p} \\
 \downarrow \Phi && \downarrow \Phi \\
 \mathbb{Z}_p^{1+3(2^{n}-1)} & \underset{\textbf{T}}{\longrightarrow } & \mathbb{Z}_p^{1+3(2^{n}-1)}
\end{array}
\end{eqnarray*}
commutes, where $\textbf{T}y = T_{f} y \mod p$ for every $y\in
\mathbb{F}^{\Gamma^{2}}_{p}.$ Since  $\Phi$ is a one-to-one
correspondence, the following statements are equivalent
\begin{enumerate}
    \item $T_{f}$ is reversible;
    \item $M_{R}^{(n)}$ is invertible over $\mathbb{F}_p$;
    \item 0 is not an eigenvalue of $M_{R}^{(n)}$ over
    $\mathbb{F}_p$;
    \item The matrix $M_{R}^{(n)}$ has a full rank.
\end{enumerate}
\end{thm}

\subsection{Illustrative Examples: Reversible}
One can compute the determinant of the rule matrix $M_R^{(n)}$ for
some random $a,b,c,d\in \mathbb{F}_{p}^{*}$ and the levels $n$ of
Cayley tree as follows:
$$
\det(M_R^{(2)})=d^4 (d^2-c (2 (a + b) + c) ) (d^2-(a + b) c)^2
$$
$$
\det(M_R^{(3)})=d^8(d^2-(a+b)c)^3(d^2-2(a+b)c)^2((a+b)c^2(a+b+c)-c(3(a+b)+c)d^2+d^4).
$$
 We have seen that the CAs are reversible for some given
values $a,b,c,d\in \mathbb{F}_{p}^{*}$ and $n$, for some values
the CAs are irreversible.

Notably, the eigenvalues of the matrix $M_R^{(2)}$ are
$d,d,d,d,d-\sqrt{a c+b c},d-\sqrt{a c+b c},\sqrt{a c+b
c}+d,\sqrt{a c+b c}+d,d-\sqrt{2 a c+2 b c+c^2},\sqrt{2 a c+2 b
c+c^2}+d$, respectively. The last situation reveals that the
necessary and sufficient conditions for the matrix $M_R^{(2)}$
being invertible are
$$
\left\{
\begin{array}{l}
 \sqrt{ac+bc}\neq d(\text{mod}p); \\
 \sqrt{2ac+2bc+c^2}\neq d(\text{mod}p).
\end{array}
\right.
$$
Therefore, the CA corresponding to the matrix $M_R^{(2)}$ is
reversible if and only if $$ \left\{
\begin{array}{l}
 \sqrt{ac+bc}\neq d(\text{mod}p); \\
 \sqrt{2ac+2bc+c^2}\neq d(\text{mod}p).
\end{array}
\right.
$$

In the Table \ref{tab:1}, we examine under what conditions these
linear transformations are invertible, and check invertibility for
a list of parameters using computations by means of "Mathematica".
For example, if we take as $a=b=c=d=1$ and $n=3$, then we can see
that the CAs are irreversible for prime numbers $p<47$. The
reversibility of finite CAs on Cayley tree of order two is
determined for some given values of $a,b,c,d\in
\mathbb{F}_{p}^{*}$ and the levels $n$ of Cayley tree. One can
fully characterize reversibility of finite cellular automata with
NBC determined by nearest neighbor rule on Cayley tree by
computing the determinant of the matrix in the Eq. \eqref{rm}.
Also, one can study the reversibility of finite CAs via rank of
the matrix in the Eq. \eqref{rm} (see \cite{SAS}).\begin{center}
\begin{table}
\caption{The reversibility of finite CAs for some given
$a,b,c,d\in \mathbb{F}_{p}^{*}$ and the levels $n=2,3$ of Cayley
tree of order two.}\label{tab:1}
\begin{tabular}{|c|c|c|c|c|c|c|}\hline
  $a$& $b$& $c$& $d$& $n$& $p$& reversibility of finite CA\\\hline
  1& 1& 1& 1& 2&2& irreversible \\\hline
  1& 1& 1& 1&2&3,5,...,101 & reversible\\\hline
 2&1& 5&2&2&17& irreversible\\\hline
 2& 1& 3& 2& 2& 17& reversible \\\hline
  2& 3& 4& 3& 2& 11& irreversible \\\hline
   1& 1& 1& 1& 3&3& irreversible \\\hline
   2& 2& 3& 3& 3&5& irreversible \\\hline
  2& 1& 1& 3& 3&5& reversible \\\hline
  2& 2& 3& 3& 3&7,11,13,19,23,29& reversible \\\hline
\end{tabular}
\end{table}
\end{center}

\section{The measure entropy of the CA on Cayley tree}
In this section we study the measure entropy of cellular automata
defined by local rules in \eqref{CA} and \eqref{CA1} on Cayley
tree of order two. In order to state our result, we first recall
necessary definitions. Let $(X, \mathcal{B},\mu,T)$ be a
measure-theoretical dynamical system. If
$\alpha=\{A_1,\ldots,A_n\}$ and $\beta=\{B_1,\ldots,B_m\}$ are two
measurable partitions of $X$, then $\alpha \vee \beta=\{A_i\cap
B_j:i=1,\ldots,n;j=1,\ldots,m\}$ is the partition of $X$. Also,
$T^{-1}\alpha$ is the partition of $X$ and
$T^{-1}\alpha=\{T^{-1}A_1,\ldots,T^{-1}A_n\}$ (see \cite{DGS,W}
for details).
\begin{defin}
Let $\alpha$ be a measurable partition of $X$. The quantity
$$
H_{\mu}(\alpha)=-{\underset{A\in \alpha}{\sum}}\mu(A)\log \mu(A)
$$
is called the entropy of the partition $\alpha$. The logarithm is
usually taken to the base 2. Let $\alpha$ be a partition with
finite entropy, then the quantity
$$
h_{\mu}(T,\alpha)=\overset{}{\underset{n\rightarrow\infty}{\lim}}
\frac{1}{n}H_{\mu}(\overset{n-1}{\underset{i=0}{\bigvee}}T^{-i}\alpha)
$$
is called the entropy of $\alpha$ with respect to $T$. The
quantity
\begin{equation}\label{ent}
h_{\mu}(T)=\overset{}{\underset{\alpha}{\sup}}\{h_{\mu}(T,\alpha):
\alpha\ \text{is a partition with}\ H_{\mu}(\alpha)<\infty\}
\end{equation} is called the measure-theoretical entropy of $(X,
\mathcal{B},\mu,T)$, the entropy of $T$ (with respect to $\mu$).
\end{defin}


Let $\pi=\{\pi_{0},\pi_{1},\ldots,\pi_{p-1}\}$ be a probability
vector. Recall that the Bernoulli measure is defined as follows:
$$
\mu_{\pi }(_0[i_0,\ldots,i_{k}]_k)=\pi_{i_{1}}\pi_{i_{0}}\ldots
\pi_{i_{k}},
$$
where $_0[i_0,\ldots,i_{k}]_k$ is a cylinder set (see \cite{DGS,W}
for details). If we take the Bernoulli measure as
$$\mu_{\pi
}(_0[i_0,\ldots,i_{k}]_k)=\frac{1}{p}\frac{1}{p}\cdots\frac{1}{p}=\frac{1}{p^{k+1}},
$$
then the measure is called uniform Bernoulli measure, i.e., for
all $i\in \mathbb{F}_{p}$, $\mu_{\pi}(_{0}[i])=\frac{1}{p}$, then
$\mu_{\pi }$ is the uniform Bernoulli measure on the space
$\mathbb{F}^{\Gamma^{2}}_{p}$. In this paper, we consider uniform
Bernoulli measure.

It is clear that due to $(a,p)=1, (b,p)=1,(c,p)=1$ and $(d,p)=1$,
the rules given in the Eqs. \eqref{CA} and \eqref{CA1} are
bipermutative. The following Theorems have been proved:

\begin{thm}\label{permut} \cite{AH} Any left-permutative (right-permutative) cellular is surjective
(see \cite{S1} for details).
\end{thm}
\begin{thm}\label{ubm}\cite{BKM} If a cellular automaton is
surjective then it preserves a uniform Bernoulli measure.
\end{thm}

D'amico \emph{et al.} \cite{DMM} have proved that for
$D$-dimensional linear CA with $D\geq 2$ the topological entropy
must be 0 or infinity (see \cite{MaMa}). In the one-dimensional
case, the measure theoretical entropy of the cellular automata is
finite \cite{A1,BKM}. In the following theorem, we prove that the
linear CA on Cayley tree of order two has infinite entropy.

Let us choose $a,b,c,d\in \mathbb{F}_{p}^{*}$ such that the
cellular automata $T_f$ defined in the Eq. \eqref{mapI} is
measure-preserving function with respect to (w.r.t.) the uniform
Bernoulli measure on the space $\mathbb{F}^{\Gamma^{2}}_{p}$. Then
we have the following theorem.

\begin{thm}
Let $T_f$ be cellular automata defined by local rules in
\eqref{CA} and \eqref{CA1} on Cayley tree of order two over the
field $\mathbb{F}_p$. Then the measure theoretical entropy of
$T_f$ w.r.t. the uniform Bernoulli measure on the space
$\mathbb{F}^{\Gamma^{2}}_{p}$ is infinity.
\end{thm}
\begin{proof}
From theorems \ref{permut} and \ref{ubm}, we note that $\mu_{\pi
}$ is a $T_f$-invariant measure. Let the zero-time partition be
given as $ \xi(0,1)= \{_{0}[0],_{0}[1],\ldots,_{0}[p-1]\},$ we put
$\xi(-i,i)=\bigvee_{u=-i}^i\sigma^{-u}\xi$, where $\sigma$ is the
shift map. Since $T_f$ is permutative, one has
\begin{eqnarray*}
\bigvee_{k=0}^{n-1}T_f^{-k}\xi(0,1)&=& \xi(0,1+3(2^{n}-1)).
\end{eqnarray*}
From the definition of measure theoretical entropy w.r.t the
measure, we get
\begin{eqnarray*}
h_{\mu_{\pi
}}(T_f,\xi(0,1))&=&{\underset{n\rightarrow\infty}{\lim}}
\frac{1}{n}H_{\mu_{\pi }}(\overset{n-1}{\underset{k=0}{\bigvee}}T_f^{-k}\xi(0,1))\\
&=&-{\underset{n\rightarrow\infty}{\lim}}
\frac{1}{n}{\underset{A\in \xi(0,1+3(2^{n}-1))}{\sum}}\mu_{\pi }(A)\log \mu_{\pi }(A) \\
&=&-{\underset{n\rightarrow\infty}{\lim}}
\frac{1}{n}p^{1+3(2^{n}-1)}\frac{1}{p^{1+3(2^{n}-1)}}\log
\frac{1}{p^{1+3(2^{n}-1)}}\\
&=&{\underset{n\rightarrow\infty}{\lim}}
\frac{1}{n}(1+3(2^{n}-1))\log p=\infty.
\end{eqnarray*}
Therefore, from the Eq. \eqref{ent}, one can conclude that
$h_{\mu_{\pi }}(T_f)=\infty.$
\end{proof}
\begin{rem} If we choose the probability vector as
$\pi=(1,0,\ldots,0)$, then\\ $h_{\mu_{\pi }}(T_f,\xi(0,1))=0$.
\end{rem}
\section{Conclusions}

In this short paper, firstly we have defined linear cellular
automata on Cayley tree of order 2. We have constructed the rule
matrix corresponding to finite cellular automata on Cayley tree by
using matrix algebra built on the field $\mathbb F_p$ (the set of
prime numbers modulo $p$). Further, we have discussed the
reversibility problem of this cellular automata. Lastly, we have
studied the measure theoretical entropy of the cellular automata
on Cayley tree.

To the best knowledge of the author, it is believed that this is
the first instance in the literature where such a connection is
established. Thus, this connection between cellular automata and
Cayley tree leads to many questions and applications that wait to
be explored.

Using the methods in the references \cite{SUAS,USAS, USAS1, SUA},
we will demonstrate  the application of linear rules on image
matrix which forms the basis of self replicating and self-similar
patterns in image processing on Cayley tree. Also, investigation
of CA on Cayley tree with more higher orders will be studied in
the future works.




%



\end{document}